# On Certain Classes of Rectangular Designs

Mithilesh Kumar Singh[1] and Shyam Saurabh[2]
[1]Ranchi University, Ranchi, India
[2]Tata College, Chaibasa, India

**Abstract**

Rectangular designs are classified as regular, Latin regular, semi–regular, Latin semi–regular and singular designs. Some series of self–dual as well as $\alpha$–resolvable designs are obtained using matrix approaches which belong to the above classes.

In every construction we obtain a matrix $N$ whose blocks are square $(0,1)$–matrices such that $N$ becomes the incidence matrix of a rectangular design. The method is the reverse of the well–known tactical decomposition of the incidence matrix of a known design. Authors have already obtained some series of Group Divisible and Latin square designs using this method. Tactical decomposable designs are of great interest because of their connections with automorphisms of designs, see Bekar *et al.* (1982). The rectangular designs constructed here are of statistical as well as combinatorial interest.

## 1. Introduction

Since the introduction of rectangular designs by Vartak (1955) workers in the design of experiments have forwarded different methods to construct them. However even a broad classification of rectangular designs have not been attempted. In this paper we identify some classes of rectangular designs and construct some of them.

*Notations:* $I_n$ denotes the identity matrix of order $n$, $J_{t \times u}$ is the $t \times u$ matrix all of whose entries are 1, in particular $J_n = J_{n \times n}$. Also $A^c = J_n - A$, $I_n^c = J_n - I_n$, $A^T$ is the transpose of a matrix $A$, $\alpha = \text{circ}(0,1,0,\ldots,0)$ is the basic circulant matrix of order $n$ such that $\alpha^n = I_n$ and $A \times B$ denotes Kronecker product of two matrices $A$ and $B$. Some relevant definitions in the context of this paper are as follows:

### 1.1 Balanced incomplete block design

A balanced incomplete block design (BIBD) or a $2 - (v, k, \lambda)$ design is an arrangement of $v$ treatments in $b$ blocks, each of size $k$ ($<v$) such that every treatment occurs in exactly $r$ blocks and any two distinct treatments occur together in $\lambda$ blocks. $v$, $b$, $r$, $k$, $\lambda$ are called parameters of the BIBD and they satisfy the relations: $bk = vr$; $r(k-1) = \lambda(v-1)$. A BIBD is symmetric (SBIBD) if $v=b$ and is self – complementary if $v = 2k$.

### 1.2 Rectangular design

Let $v = mn$ treatments be arranged in an $m$ x $n$ array $A$. A rectangular design (RD) is an arrangement of the $v = mn$ treatments in $b$ blocks each of size $k$ such that:

Corresponding Author: Shyam Saurabh
Email: shyamsaurabh785@gmail.com

1. Every treatment occurs at most once in a block;
2. Every treatment occurs in exactly *r* blocks;
3. Every pair of treatments, which are in the same row of *m* x *n* array occur together in $\lambda_1$ blocks; which are in same column occur together in $\lambda_2$ blocks; while every other pair of treatments occur together in $\lambda_3$ blocks. The treatments in the same row are first associates, the treatments in the same column are second associates and the treatments not in the same row or column are third associates.

$v = mn$, $b$, $r$, $k$, $\lambda_1$, $\lambda_2$ and $\lambda_3$ are known as parameters of the RD and they satisfy the relations: $bk = vr; (n-1)\lambda_1 + (m-1)\lambda_2 + (n-1)(m-1)\lambda_3 = r(k-1)$. (1)

Let *N* be the incidence matrix of the RD then the structure of $NN^T$ is given as:

$$NN^T = \begin{pmatrix} (r-\lambda_1)I_n + \lambda_1 J_n & (\lambda_2-\lambda_3)I_n + \lambda_3 J_n & \cdots & (\lambda_2-\lambda_3)I_n + \lambda_3 J_n \\ (\lambda_2-\lambda_3)I_n + \lambda_3 J_n & (r-\lambda_1)I_n + \lambda_1 J_n & \cdots & (\lambda_2-\lambda_3)I_n + \lambda_3 J_n \\ \vdots & \vdots & \ddots & \vdots \\ (\lambda_2-\lambda_3)I_n + \lambda_3 J_n & (\lambda_2-\lambda_3)I_n + \lambda_3 J_n & \cdots & (r-\lambda_1)I_n + \lambda_1 J_n \end{pmatrix}$$

$$= r(I_m \times I_n) + \lambda_1(I_m \times I_n^c) + \lambda_2(I_m^c \times I_n) + \lambda_3(I_m^c \times I_n^c).$$

The *m* x *n* array *A* is given as:
$$\begin{matrix} 1 & 2 & 3 & \cdots & n \\ n+1 & n+2 & n+3 & \cdots & 2n \\ \vdots & \vdots & \vdots & \vdots & \vdots \\ (m-1)n+1 & (m-1)n+2 & (m-1)n+3 & \cdots & mn \end{matrix}.$$

An RD with $v = mn$ treatments reduces to

(i) a *group divisible design* (*GDD*) if $\lambda_1 = \lambda_3$ or $\lambda_2 = \lambda_3$. Furthermore, if $r-\lambda_1 = 0$ then the GDD is singular (*S*); if $r - \lambda_1 > 0$ and $rk - v\lambda_2 = 0$ then it is semi–regular (*SR*); and if $r - \lambda_1 > 0$ and $rk - v\lambda_2 > 0$, the design is regular (*R*).

(ii) an $L_2$ – *type Latin square design* ($L_2$ –*type design*) if $m = n, \lambda_1 = \lambda_2$, see Nigam *et al.* (1988, p. 58).

For details on these designs, see Raghavarao (1971), Raghavarao and Padgett (2005) and Dey (1986, 2010), see also Table 1 of Section 2.

### 1.3 α − Resolvable design

A block design *D* (*v, b, r, k*) whose *b* blocks can be divided into $t = r/\alpha$ classes, each of size $\beta = v\alpha/k$ and such that in each class of $\beta$ blocks every element of *D* is replicated $\alpha$ times, is called an $\alpha$ − resolvable design. The classes are called parallel or resolution classes. If $\alpha = 1$ then the design is said to be resolvable. Further a resolvable block design is said to be affine resolvable if any two distinct blocks from the same resolution class intersect in the same number, say $q_1$, of treatments; whereas any two blocks belonging to different resolution classes intersect in the same number, say $q_2$, of treatments.

Alternatively, if the incidence matrix $N$ of a block design $D$ ($v$, $b$, $r$, $k$) may be partitioned in to submatrices as: $N = (N_1|N_2|\cdots|N_t)$ where each $N_i (1 \leq i \leq t)$ is a $v \times v\alpha/k$ matrix such that each row sum of $N_i$ is $\alpha$ then the design is $\alpha -$ resolvable.

## 1.4 Tactical decomposable designs

Let (0, 1) – matrix $N$ have a decomposition $N = [N_{ij}]_{\substack{i=1,2,\ldots,s \\ j=1,2,\ldots,t}}$ where $N_{ij}$ are submatrices of $N$ of suitable sizes. The decomposition is called row tactical if row sum of $N_{ij}$ is $r_{ij}$ and column tactical if the column sum of $N_{ij}$ is $k_{ij}$; $i = 1,2,3,\ldots,s; j = 1,2,3,\ldots,t$ and tactical if it is row as well as column tactical. If $N$ is the incidence matrix of a block design $D$, $D$ is called the row (column) tactical decomposable. $D$ is called uniform row (column) tactical decomposable if $r_{ij} = \alpha (k_{ij} = \beta) \forall i,j$. If each $N_{ij}$ is an $n \times n$ matrix, $D$ is called square tactical decomposable design, STD ($n$).

## 1.5 Global decomposable Self – dual RD (SDRD)

The RD whose incidence matrix is $N^T$ is called the dual of the design with incidence matrix $N$. If the design with incidence matrix $N^T$ has the same parameters as the design with incidence matrix $N$, the original design is self–dual. If the two designs are also isomorphic, we say that the original design is truly self–dual.

Let $N = [N_{ij}]_{i,j=1,2,\ldots,m}$ be the incidence matrix of an SDRD D. Let $R_i$ be $i^{\text{th}}$ row of submatrices $N_{ij} (j = 1,2,\ldots,m)$. $D$ will be called global decomposable if each $N_{ij}$ is an SBIBD such that

$(i) R_i^2 = R_i R_i^T = \sum_{k=1}^{m} N_{ik} N_{ik}^T = kI_n + \lambda_1 I_n^c; i = 1,2,3,\ldots,m;$

$(ii) R_i \cdot R_j = R_i R_j^T = \sum_{k=1}^{m} N_{ik} N_{jk}^T = \lambda_2 I_n + \lambda_3 I_n^c; i \neq j.$

An SDRD $D$ has 3 (possibly 2) intersection numbers and so it may be called nearly affine design [See Ionin and Shrikhande (1995)]. In fact if $D$ has parameters $v, k, \lambda_1, \lambda_2, \lambda_3, m, n$ then the blocks of $D$ can be arranged in an $m \times n$ array such that any two distinct blocks in the same row meet in $\lambda_1$ treatments, those in the same column meet in $\lambda_2$ treatments and other pair of distinct blocks meet in $\lambda_3$ treatments. For details on global decomposable designs, see Ionin and Shrikhande (2006).

## 1.6 Skew – Hadamard design

Following Koukouvinos and Stylianou (2008), a (1,–1) – matrix $H$ of order $4t$ is called a skew Hadamard matrix (skew H–matrix) if

$(i) HH^T = 4tI_{4t}$ $(ii) H - I_{4t}$ is a skew symmetric matrix.

A skew $H$ – matrix is in normalized form if its first row as well as first column contains entirely 1's. The matrix of order $4t$–1 obtained by deleting $1^{\text{st}}$ row and $1^{\text{st}}$ column of a normalized skew

H–matrix $H$ is called the core $C$ of $H$. Then $N = (J + C)/2$ is the incidence matrix of a $(4t–1, 2t–1, t–1)$ – design called skew– Hadamard design and it satisfies the following relations:

$$N + N^T = I^c_{4t-1}, NN^T = N^TN = (2t - 1)I_{4t-1} + (t - 1)I^c_{4t-1};$$

$$N^2 + N = (N^T)^2 + N^T = tI^c_{4t-1}, N^2 + (N^T)^2 = (2t - 1)I^c_{4t-1}.$$

## 2. Classification of RDs

Let $N$ be the incidence matrix of a RD $D$ with parameters: $v, r, k, b, \lambda_1, \lambda_2, \lambda_3, m, n$. Assume that $\theta_1, \theta_2, \theta_3$ are the eigenvalues of $NN^T$ given by

$$\theta_1 = r - \lambda_1 + (m - 1)(\lambda_2 - \lambda_3), \theta_2 = r - \lambda_2 + (n - 1)(\lambda_1 - \lambda_3), \theta_3 = r - \lambda_1 - \lambda_2 + \lambda_3 \quad (2)$$

with multiplicities $n - 1, m - 1, (m - 1)(n - 1)$ respectively. (3)

Then $D$ will be called

$(r_1)$ Regular (R) if $\theta_1, \theta_2, \theta_3 > 0$;
$(r_2)$ Semi – regular of type I [SR(I)] if $\theta_1 = 0; \theta_2, \theta_3 > 0$;
$(r_3)$ Semi–regular of type II [SR(II)] if $\theta_2 = 0; \theta_1, \theta_3 > 0$;
$(r_4)$ Latin semi – regular [(LSR)] if $\theta_1 = \theta_2 = 0, \theta_3 > 0$;
$(r_5)$ Singular [(S)] if $\theta_3 = 0$.

Clearly the above five classes of RDs are mutually exclusive and exhaustive. $(r_1)$ has a special subtype $(r'_1)$: RD is Latin regular (LR) if $\theta_1 = \theta_2(> 0)$.

The following proposition provides the basis for the above classification for RD D:

**Proposition 1:** (a) If $\theta_1 = \theta_3(> 0), \theta_2 > 0$ or $\theta_2 = \theta_3(> 0), \theta_1 > 0$ then $(r_1)$ reduces to a regular GD.
(b) If $\theta_1 = \theta_2(> 0), \theta_3 > 0$ and $m = n$ then $(r_1)$ or $(r'_1)$ reduces to an $L_2$ −type design.
(c) If $\theta_1 = 0; \theta_2 = \theta_3(> 0)$ then $(r_2)$ reduces to a semi– regular GD with $\lambda_2 = \lambda_3$.
(d) If $\theta_2 = 0; \theta_1 = \theta_3(> 0)$ then $(r_3)$ reduces to a semi– regular GD with $\lambda_1 = \lambda_3$.
(e) If $\theta_1 = \theta_2 = 0, \theta_3 > 0$ then $(r_4)$ does not reduce to a GDD in any case. However for $m = n$, it reduces to an $L_2$ − type design. Hence this RD is given the name LSR RD.
(f) If $\theta_1 = 0$ or $\theta_2 = 0$, then $(r_5)$ reduces to a singular GD.
(g) If $\theta_1 = \theta_2 = \theta_3 = 0$, then $v = k$ and hence $D$ is the complete design (a trivial singular GDD). Further if $\theta_1, \theta_2(> 0)$ and $\theta_3 = 0$, the following example shows that there is a singular RD which is not a GDD:

**Example 1:** Consider $v = n^2$ treatments arranged in an $n \times n$ array $A$ and take rows and columns of the array $A$ as blocks. Then we obtain a connected singular RD with parameters: $v = n^2, b = 2n, r = 2, k = n, \lambda_1 = \lambda_2 = 1, \lambda_3 = 0, m = n, n$ which is not a GDD.

**Remark 1:** The transpose of $m$ x $n$ array $A$ makes the interchanges: $\lambda_1 \leftrightarrow \lambda_2, m \leftrightarrow n, \theta_1 \leftrightarrow \theta_2$ and SR (I) RD↔SR (II) RD. Hence following the convention in GDD, we often will not make distinction between SR (I) and SR (II) RDs and each will be called SR RD.

The following Table 1 reveals the relation of RD with its subdesigns:

| Table 1: RDs and its Subdesigns | | | |
|---|---|---|---|
| Design | Regular RD | | Semi–regular RD |
| Subdeisgns | (i) Latin regular (LR) if $\theta_1 = \theta_2 (> 0), \theta_3 > 0$ <br> (ii) RGDD if $\theta_1 = \theta_3 (> 0), \theta_2 > 0$ or $\theta_2 = \theta_3 (> 0), \theta_1 > 0$ <br> (iii) BIBD of order $mn$ if $\theta_1 = \theta_2 = \theta_3 (> 0)$ | | SRGDD if $\theta_1 = 0; \theta_2 = \theta_3 (> 0)$ <br> or $\theta_2 = 0; \theta_1 = \theta_3 (> 0)$ |
| Design | Latin semi–regular RD | | Singular RD |
| Subdeisgns | $L_2$ –type design if $\theta_1 = \theta_2 = 0, m = n$ | | SGDD if $\theta_1 = \theta_3 = 0$ or $\theta_2 = \theta_3 = 0$ |

**Proposition 2:** Let D be an RD with parameters: $v, r, k, b, \lambda_1, \lambda_2, \lambda_3, m, n$. Then
(a) If $D$ is SR, then at least one of $\lambda_3 - \lambda_1, \lambda_3 - \lambda_2$ is positive.
(b) If $D$ is LR and one of $\lambda_3 - \lambda_1, \lambda_3 - \lambda_2$ is positive, then the other is also positive.
(c) If $D$ is LSR, then both $\lambda_3 - \lambda_1, \lambda_3 - \lambda_2$ are positive.
(d) If $D$ is S, none of $\lambda_3 - \lambda_1, \lambda_3 - \lambda_2$ is positive.

**Proof:** Proof of (a) and (d) is an easy exercise.

(b) For an LR RD: $\theta_1 = \theta_2 (> 0) \Rightarrow m(\lambda_3 - \lambda_2) = n(\lambda_3 - \lambda_1)$.

(c) For an LSR RD: $\theta_1 = \theta_2 = 0 \Rightarrow r = \lambda_1 + (m-1)(\lambda_3 - \lambda_2) = \lambda_2 + (n-1)(\lambda_3 - \lambda_1)$ (4)

Using relations (1) and (4), we obtain

$\lambda_3 - \lambda_2 = \frac{nr}{d}(mn - k) > 0; \lambda_3 - \lambda_1 = \frac{mr}{d}(mn - k) > 0$ where $d = mn(m-1)(n-1)$.

**Theorem 1:** Let $D$ be a semi – regular (II) RD with parameters: $v, r, k, b, \lambda_1, \lambda_2, \lambda_3, m, n$ and $A$ be the $m$ x $n$ array as given in definition 1.2, then
(i) $m$ divides $k$ and if $k/m = \alpha$, then $\alpha$ is the number of treatments common to any block of $D$ and any row of $n$ treatments in the $m$ x $n$ array $A$;
(ii) $D$ is uniform column tactical decomposable;
(iii) The dual of $D$ is $\alpha$ – resolvable;
(iv) If $D$ is self–dual, it is $\alpha$ –resolvable and STD $(n)$;
(v) If $D$ is self–dual, there exist $(m - 1)$ semi–regular $\alpha$–resolvable RDs $D'$ with parameters: $v' = (m - s)n, m' = m - s, n' = n, k' = \alpha(m - s), \lambda_1, \lambda_2, \lambda_3; 1 \leq s \leq m - 1$.

**Proof:** (i) Using (1) in (2), we obtain $\theta_2 = rk - m[\lambda_2 + (n-1)\lambda_3]$ (5)
Let $P_i, i = 1, 2, \ldots, m$ be the sets of points which form the rows of the array A, taken in order. Let $B_j, j = 1, 2, \ldots, b$ be the blocks of $D$. Consider a matrix: $M = [m_{ij}]_{\substack{i=1,2,\ldots,m \\ j=1,2,\ldots,b}}$, where $m_{ij} = |P_i \cap B_j|$. By simple combinatorial arguments, it follows that
$\sum_{j=1}^{b} m_{ij} = nr \; \forall i$ and $\sum_{j=1}^{b} m_{ij}(m_{ij} - 1) = n(n-1)\lambda_1$. (6)
Adding these equations, we get: $\sum_{j=1}^{b} m_{ij}^2 = nr + n(n-1)\lambda_1$ (7)

Next we have $\sum_{i=1}^{m}\sum_{j=1}^{b}(m_{ij} - \frac{k}{m})^2 = \frac{n}{m}[m(r + (n-1)\lambda_1) - kr]$ (using (6), (7) and $b = \frac{mnr}{k}$)
$= 0$ [putting $\theta_2 = 0$ in (5) and using (1)]. Hence $m_{ij} = \frac{k}{m}$.

(ii) We order the rows of the incidence matrix $N$ of $D$ such that first set of $n$ rows correspond to $P_1$, second set of $n$ rows correspond to $P_2$ and so on. This partitions $N$ into $m$ $n \times b$ matrices. By (i), the column sum of each matrix is $\frac{k}{m} = \alpha$ (say). Hence (ii) follows.

Further the dual of $D$ has the incidence matrix $N^T$ partitioned into $m$ $b \times n$ matrices each with row sum $\alpha$. Hence (iii) and (iv) follows.

(v) The incidence matrix of $D'$ can be obtained from that of $D$ by deleting any $s$ $n \times b$ matrices constructed in (ii).

**Remark 2:** If $D$ is semi–regular (I) and $n$ divides $k$, then Theorem also holds, which can be observed by interchanging $m$, $n$ as well as $\lambda_1, \lambda_2$.

**Theorem 2: Properties of a regular RD D**
(i) If $D$ is symmetric, $\theta_1^{n-1}\theta_2^{m-1}\theta_3^{(m-1)(n-1)}$ is a perfect square;
(ii) For a regular resolvable RD $b \geq v + r - 1$;
(iii) A symmetric regular RD cannot be resolvable unless $r = 1$.

**Proof:** (i) This follows from the relation: $\det(MM^T) = (detM)^2 =$ the product of eigenvalues of $MM^T = \theta_1^{n-1}\theta_2^{m-1}\theta_3^{(m-1)(n-1)}$.

The proofs of (ii) and (iii) are the same as those for a regular GDD as given in Bose and Connor (1952).

## 3. Earlier Constructions

Vartak (1955), Raghavarao and Aggarwal (1974), Kageyama and Tanaka (1981), Banerjee *et al.* (1985), Bhagwandas *et al.* (1985), Suen (1989), Sinha (1991), Kageyama and Miao (1995), Sinha and Mitra (1999), Sinha *et al.* (1993, 1996, 2002*a*), Kageyama and Sinha (2003), Bagchi (1994, 2004) have constructed RDs using various approaches. However their constructions methods may be summarized as:

| | |
|---|---|
| Vartak (1955): | From Kronecker product of $2(v, k, \lambda) -$ designs |
| Raghavaro and Aggarwal (1974): | From difference sets |
| Kageyama and Tanaka (1981): | From $2 - (v, k, \lambda)$ designs and Skew Hadamard designs |
| Bhagwandas *et al.* (1985): | From $2 - (v, k, \lambda)$ designs |
| Banerjee *et al.* (1985): | From Rectangular association scheme |
| Suen (1989) and Bagchi (1994): | From difference sets |
| Sinha (1991): | From Bhaskar Rao designs |
| Kageyama and Miao (1995): | From $2 - (v, k, \lambda)$ designs, Self–complementary $2 - (v, k, \lambda)$ designs, Normalised Hadamard matrices, Breaking up blocks |

| Sinha *et al.* (1996): | From balanced bipartite weighing designs |
| Sinha and Mitra (1999): | From nested and resolvable $2-(v,k,\lambda)$ designs |
| Sinha *et al.* (1993, 2002a): | From $2-(v,k,\lambda)$ designs, nested and self–complementary $2-(v,k,\lambda)$ designs; Semi–regular group divisible designs; Mutually orthogonal Latin squares |
| Kageyama and Sinha (2003): | From $2-(v,k,\lambda)$ design and nested $2-(v,k,\lambda)$ designs |
| Bagchi (2004): | From resolvable and almost resolvable $2-(v,k,\lambda)$ designs |
| Parihar *et al.* (2009): | From $2-(v,k,\lambda)$ designs and Hadamard matrices |

The constructions of RDs by above authors may be classified in to four classes as discussed above. Sinha *et al.* (2002) obtained some series of cyclic RDs and some examples of partial cyclic RDs may be found in Sinha and Mitra (1999). E– optimality of these designs have been discussed by Bagchi (1994, 2004) and Sinha *et al.* (2002a), among others.

## 4. Present Constructions

**Method I: From $2-(v,k,\lambda)$ and SH – designs**

**Lemma 1:** Let $N_i (i=1,2)$ be the incidence matrices of SBIBDs. Then $N_1 \times N_2$ is the incidence matrix of a global decomposable SDRD.

**Theorem 3:** The existence of a $(v_1, k_1, \lambda'_1)$ – design and a $(4t–1, 2t–1, t–1)$ – SH design implies the existence of a global decomposable STD $(4t–1)$ RD with parameters:

$$v = (4t-1)v_1, b = (4t-1)b_1, r = 2r_1(t-1) + b_1, k = 2k_1(t-1) + v_1, \lambda_1 = (t-1)r_1,$$

$$\lambda_2 = b_1 - 2r_1 + 2t\lambda'_1, \lambda_3 = r_1 + \lambda'_1(t-2), m = v_1, n = 4t-1. \qquad (8)$$

**Proof:** Let $N_1$ be the incidence matrix of a $(v_1, k_1, \lambda'_1)$ – design and $N_2$ be the incidence matrix of a $(4t–1, 2t–1, t–1)$ – SH design. Then $N_2 N_2^T = N_2^T N_2 = tI_{4t-1} + (t-1)J_{4t-1}$, $N_2 + N_2^T = (J-I)_{4t-1}$. Using these relations one can verify that $N = N_1 \times N_2 + (J_{V_1} - N_1) \times I_{4t-1}$ is the incidence matrix of the RD with parameters (8).

**Remark 3:** The Theorem 3 may be found in Kageyama and Tanaka (1981). Here we have used algebraic approach in the proof and the series is also identified as global decomposable.

**Theorem 4:** The existence of symmetric $2-(n, k_i, \lambda_i)$ designs $(i=1,2)$ whose incidence matrices $N_i, i=1,2$ satisfy the relation: $N_1 N_2^T + N_2 N_1^T = \mu_1 I_n + \mu_2 I_n^c$, implies the existence of a global decomposable STD $(n)$ SDRD with parameters: $v = b = mn, r = k = k_1 + (m-2)k_2, \lambda'_1 = \lambda_1 + (m-1)\lambda_2, \lambda'_2 = (m-2)k_2 + \mu_1, \lambda'_3 = (m-2)\lambda_2 + \mu_2, m, n.$ (9)

**Proof:** Let $N = I_m \times N_1 + I_m^c \times N_2$. Then it can be verified that:

$NN^T = N^TN = \{k_1 + (m-2)k_2\}(I_m \times I_n) + \{\lambda_1 + (m-1)\lambda_2\}(I_m \times I_n^c) + \{(m-2)k_2 + \mu_1\}(I_m^c \times I_n) + \{(m-2)\lambda_2 + \mu_2\}.$

Hence $N$ represents the incidence matrix of a SDRD with parameters (9).

In the following series (Corollaries 1 – 7) of RDs, $w = circ.(0\ 1\ 0)$ and $N_{4t-1} = circ.(0\ 1\ 1\ 0\ 1\ 0\ 0)$ denote incidence matrices of SH–designs with $v = 3$ and $v = 4t - 1$ treatments respectively.

**Corollary 1:** For $N_1 = N_{4t-1}, N_2 = I_{4t-1}$, we get a series of regular (for $t > 1$) SDRDs with parameters: $v = m(4t - 1), k = m + 2t - 2, \lambda_1 = t - 1, \lambda_2 = m - 2, \lambda_3 = 1, m, n = 4t - 1$.

**Remark 4:** When $4t - 1$ is a prime or prime power, the series in Corollary 1 has the same parameters as in the series given by Sinha *et al.* (2002a) in Theorem 2.4.

**Corollary 2:** For $N_1 = I_{4t-1} + N_{4t-1}, N_2 = I_{4t-1}$, we get a series of regular SDRDs with parameters: $v = m(4t - 1), k = m + 2t - 1, \lambda_1 = t, \lambda_2 = m, \lambda_3 = 1, m, n = 4t - 1$.

**Corollary 3:** For $N_1 = I_{4t-1}, N_2 = N_{4t-1}$, we get a series of SDRDs with parameters: $v = m(4t-1), k = (m-1)(2t-1) + 1, \lambda_1 = (m-1)(t-1), \lambda_2 = (m-2)(2t-1), \lambda_3 = (m-2)(t-1) + 1, m, n = 4t - 1$.

**Remark 5:** The above RD is regular when $t>1$ and semi–regular when $t = 1$. When $4t - 1$ is a prime or prime power, the parameters of the series in Corollary 3 contain those of the series in Theorems 2.1, 2.2, 2.3 in Sinha *et al.* (2002a) and the series of the Theorem 2.4 given in Sinha *et al.* (1993).

**Corollary 4:** For $N_1 = I_m, N_2 = w$, we get the series of semi– regular resolvable SDRDs with parameters: $v = 3m, k = m, \lambda_1 = 0, \lambda_2 = m - 2, \lambda_3 = 1, m, n = 3$.

**Corollary 5:** For $N_1 = N_{4t-1}, N_2 = N_{4t-1}^T$, we obtain a series of $(2t - 1)$ – resolvable semi–regular (II) SDRD with parameters: $v = m(4t - 1), k = m(2t - 1), \lambda_1 = m(t - 1), \lambda_2 = (m - 2)(2t - 1), \lambda_3 = m(t - 1) + 1, m, n = 4t - 1$.

**Corollary 6:** For $N_1 = I_{4t-1} + N_{4t-1}, N_2 = I_{4t-1} + N_{4t-1}^T$, the series is 2t–resolvable semi–regular SDRD with parameters: $v = m(4t - 1), k = 2mt, \lambda_1 = mt, \lambda_2 = 2t(m - 2) + 2, \lambda_3 = mt + 1, m, n = 4t - 1$.

**Corollary 7:** When $N_1 = I_{4t-1} + N_{4t-1}, N_2 = N_{4t-1}^T$, we get a series of SDRD with parameters: $v = m(4t - 1), k = m(2t - 1) + 1, \lambda_1 = m(t - 1) + 1, \lambda_2 = (m - 2)(2t - 1), \lambda_3 = m(t - 1) + 2, m, n = 4t - 1$.

**Theorem 5:** The existence of a $2 - (v, k, \lambda)$ design satisfying $(n - 1)^2 = \frac{kr}{r-\lambda}$ (10)

implies the existence of an LR RD with parameters: $v^* = vn, b^* = bn, r^* = r(n-1), k^* = k(n-1), \lambda_1^* = (n-2)r, \lambda_2^* = \lambda(n-1), \lambda_3^* = \lambda(n-2), m = v, n.$ (11)

**Proof:** Let $M$ be the incidence matrix of a $2 - (v, k, \lambda)$ design satisfying (10). Using the relation: $m(\lambda_3^* - \lambda_2^*) = n(\lambda_3^* - \lambda_1^*)$ for an LR RD [vide Proposition 2(c)]; it is easy to verify that $N = M \times I_n^c$ is the incidence matrix of the LR RD with parameters (11).

**Remark 6:** The $2 - (v, k, \lambda)$ designs given in Mathon and Rosa (2007): MR2, MR13, MR35, MR36, MR37, MR66, MR76, MR84, MR142, MR145 and MR166 satisfy the condition of Theorem 5 and hence are suitable for the construction of LR RDs. Here MR$X$ denotes the designs listed in Mathon and Rosa (2007) with number $X$. However we obtain only one LR RD under the range $r, k \leq 10$ using MR2: $2 - (9,3,1)$ design in Theorem 5 whose parameters are: $v = 27, b = 36, r = 8, k = 6, \lambda_1 = 4, \lambda_2 = 2, \lambda_3 = 1, m = 9, n = 3$. This is the design number 28 of Kageyama and Sinha (2003).

**Method II: Using Circulant Matrices**

**Theorem 6:** There exists a STD (4t–1) SDRD with parameters: $v = 5(4t - 1), k = 8t - 3, \lambda_1 = 4(t - 1), \lambda_2 = 2t - 1, \lambda_3 = 3t - 1, m = 5, n = 4t - 1$. (12)

**Proof:** Consider the ordered set $\{0, 1, 2, 3, 4\}$ which has the elements of the finite field GF(5). Replacing the squares by 1 and non–squares by –1 in GF(5) we obtain $\{0\ 1\ -1\ -1\ 1\}$ which we take as the first row of a circulant matrix C. i.e. C = circ. $(0, 1, -1, -1, 1)$. Further replacing 0 by $I_{4t-1}$, 1 by $N$ and –1 by $N^T$ in C where $N$ is the incidence matrix of a $(4t–1, 2t–1, t–1)$ – SH design we obtain

$$N_1 = \begin{pmatrix} I & N & N & N^T & N^T \\ N & I & N^T & N^T & N \\ N & N^T & I & N & N^T \\ N^T & N^T & N & I & N \\ N^T & N & N^T & N & I \end{pmatrix}$$ which represents the incidence matrix of STD(4t–1) SDRD with parameters (12).

**Remark 7:** $N_1^c$ represents a SDRD with parameters: $v = 5(4t - 1), k = 2(6t - 1), \lambda_1 = 8t - 3, \lambda_2 = 6t, \lambda_3 = 7t, m = 5, n = 4t - 1$.

The above two series of RDs are regular when $t > 1$ and semi–regular when $t = 1$. Some RDs from ciculant matrices are given below:

**(1)** If $A_1 = \alpha + \alpha^6, A_2 = \alpha^2 + \alpha^5, A_3 = \alpha^3 + \alpha^4$ where $\alpha$ is a basic circulant matrix of order 7 then

**Example 2:** $N = circ\ (A_1, A_2, A_3)$ is the incidence matrix of a semi–regular 2–resolvable RD with parameters: $v = b = 21, r = k = 6, \lambda_1 = 0, \lambda_2 = 1, \lambda_3 = 2, m = 7, n = 3$.

**Example 3:** $N = circ\ (I_7 + A_1, I_7 + A_2, I_7 + A_3)$ is the incidence matrix of semi–regular 3–resolvable SDRD with parameters: $v = b = 21, r = k = 9, \lambda_1 = 3, \lambda_2 = 3, \lambda_3 = 4, m = 3, n = 7$.

**Example 4:** Let $\alpha$ be the basic circulant matrix of order 11. If $A_i = \alpha^i + \alpha^{11-i}, i = 1, 2, \ldots, 5$ then $N = \begin{bmatrix} A_1 & A_2 & A_3 & A_4 & A_5 \\ A_2 & A_4 & A_5 & A_3 & A_1 \end{bmatrix}$ is the incidence matrix of a semi–regular 2–resolvable RD with parameters: $v = 22, b = 55, r = 10, k = 4, \lambda_1 = 0, \lambda_2 = 1, \lambda_3 = 2, m = 11, n = 2$.

**Method III: From Mutually orthogonal Latin squares**

A Latin square of order $n$ is an $n \times n$ array on $n$ symbols such that each of the $n$ synbols occurs exactly once in each row and each column. The join of two Latin squares $A = [a_{ij}]_{1 \leq i, j \leq n}$ and $B = [b_{ij}]_{1 \leq i, j \leq n}$ is the $n \times n$ array whose $(i, j)$ – th entry is the ordered pair $(a_{ij}, b_{ij})$. Two Latin squares are orthogonal if the join of A and B contains every ordered pair exactly once. A set of Latin squares are mutually orthogonal (MOLS) if they are pairwise orthogonal.

**Theorem 7:** The existence of $m$ MOLS of side $n$ implies the existence of a resolvable semi regular STD $(n)$ RD with parameters: $v = mn, b = n(n-1), r = n-1, k = m, \lambda_1 = \lambda_2 = 0, \lambda_3 = 1, m, n.$ (13)

**Proof:** Have $m$ MOLS of order $n$ each with the first row $(1, 2, \ldots, n)$. Deleting the first row we get $m$ rectangular MOLS $L_i (i = 1, 2, \ldots, m)$. From these MOLS we obtain STD $(n)$ matrix as follows:

Let $L_i = \begin{bmatrix} a_{21} & a_{22} & \cdots & a_{2n} \\ \vdots & \vdots & \cdots & \vdots \\ a_{n1} & a_{n2} & \cdots & a_{nn} \end{bmatrix}$. From the first row obtain an $n \times n$ permutation matrix $N_{i2}$ whose $i^{\text{th}}$ column has 1 at $a_{2i}$ – th $(i = 1, 2, \ldots, n)$ position and 0 elsewhere. Repeating the process for other rows of $L_i$, we get a sequence of permutation matrices $N_{ij}; j = 2, 3, \ldots, n$ which we arrange in a row. The $m$ rows corresponding to $L_i (i = 1, 2, \ldots, m)$ constitute a block matrix $N = [N_{ij}]_{\substack{i=1,2,\ldots,m \\ j=2,3,\ldots,n}}$. We claim that $N$ is the required STD $(n)$ matrix. It is sufficient to show that

(i) $\sum_{j=2}^{n} N_{ij} N_{ij}^T = (n-1) I_n$  (ii) $\sum_{j=2}^{n} N_{ij} N_{kj}^T = I_n^c$.

Since each $N_{ij}$ is a permutation matrix $N_{ij} N_{ij}^T = I_n$, (i) follows. Next we note that If $N_{ij}$ and $N_{kj}$ are the permutation matrices corresponding to the $j$th row $(a_{i1}, a_{i2}, \ldots, a_{in})$ of $L_i$ and $j$th row $(b_{k1}, b_{k2}, \ldots, b_{kn})$ of $L_k$ respectively, then $N_{ij} N_{kj}^T$ is the adjacency matrix of the binary relation $\{(a_{i1}, b_{k1}), (a_{i2}, b_{k2}) \ldots (a_{in}, b_{kn})\}$. Since $L_i$ and $L_k$ are mutually orthogonal rectangular moles each with the only first row $(1, 2, \ldots, n)$ missing, $\sum_{j=2}^{n} N_{ij} N_{kj}^T$ is the adjacency matrix of the binary relation $\cup \{(a_{i1}, b_{k1}), (a_{i2}, b_{k2}) \ldots (a_{in}, b_{kn})\} = \{(1,1), (2,2), \ldots, (n,n)\}^c$. Hence $\sum_{j=2}^{n} N_{ij} N_{kj}^T = I_n^c$ and consequently $N$ is the incidence matrix of an RD with parameters (13). Since row sum of each $N_{ij}$ is 1, the RD is resolvable.

**Corollary 8:** There exists a semi–regular STD $(q)$ resolvable SDRD with parameters: $v = q^2 - q, k = q - 1, \lambda_1 = \lambda_2 = 0, \lambda_3 = 1, m = q - 1, n = q$ where $q$ is a prime power.

Corollary 8 follows from the fact that there are $(q - 1)$ MOLS of order $q$.

**Corollary 9:** There exists a semi–regular 2 – resolvable STD (q) SDRD with parameters: $v = q^2 - q, k = 2(q-1), \lambda_1 = 2, \lambda_2 = q - 1, \lambda_3 = 3, m = q - 1, n = q$ (14) where $q$ is a prime power.

**Proof:** Let $N = [N_{ij}]_{\substack{i=1,2,\dots,m \\ j=1,2,\dots,n}}$ be the incidence matrix of RDs in Corollary (8) then it can be verified that $M = [I_q + N_{ij}]_{\substack{i=1,2,\dots,m \\ j=1,2,\dots,n}}$ is the incidence matrix of the RD with parameters (14).

**Remark 8:** Theorem 7 may be found in Sinha *et al.* (2002a) when $n$ is a prime or prime power. Here the proof is given using algebraic approach which gives tactical decomposition of the designs. Further by the above process as given in Theorem 7, we can obtain non symmetric resolvable STD RD when $v$ is not a prime or prime power.

**Example 5:** For $v = 10$, we have only two MOLS. Hence a resolvable RD with parameters: $v = 20, b = 90, r = 9, k = 2, \lambda_1 = \lambda_2 = 0, \lambda_3 = 1, m = 2, n = 10$ can be constructed.

**Corollary 10:** There exists a semi–regular $(q - 1)$ – resolvable STD $(q)$ SDRD with parameters: $v = q^2 - q, k = (q-1)^2, \lambda_1 = \lambda_2 = (q-1)(q-2), \lambda_3 = (q-1)(q-2) + 1, m = q - 1, n = q$ where $q$ is a prime power.

Corollary 10 is the complementary of the RD given in Corollary 8.

If $q=p$ is a prime, the alternative proof of Corollaries 8 and 9 are as follows:

**Corollary 11:** Let $p$ be a prime number and let $C$ be the basic circulant matrix of order $p$. Then

(i) $N = [C^{ij}]_{i,j=1,2,\dots,p-1}$ of a STD ($p$) resolvable semi– regular SDRD with parameters: $v = p^2 - p, k = p - 1, \lambda_1 = \lambda_2 = 0, \lambda_3 = 1, m = p - 1, n = p$ where $C^{ij} = (C^i)^j$.

(ii) $N = [I_p + C^{ij}]_{i,j=1,2,\dots,p-1}$ of a 2 – resolvable semi– regular STD ($p$) SDRD with parameters: $v = p^2 - p, k = 2(p-1), \lambda_1 = 2, \lambda_2 = p - 1, \lambda_3 = 3, m = p - 1, n = p$ where $C^{ij} = (C^i)^j$.

**Method IV: From Strongly regular graphs**

A *strongly regular graph* with parameters $(v, k, \lambda, \mu)$, $srg(v, k, \lambda, \mu)$ is a finite graph on $v$ vertices, without loops or multiple edges, regular of degree $k$ (with $0 < k < v - 1$, so that there are both edges and nonedges), and such that any two distinct vertices have $\lambda$ common neighbors when they are adjacent, and $\mu$ common neighbors when they are nonadjacent.

**Theorem 8:** The existence of $srg(v, k, \lambda, \mu)$ satisfying $\mu = \lambda + 1$ implies the existence of the STD ($v$) SDRDs with parameters:

(i) $v_1 = 2v, k_1 = v - 1, \lambda_1 = v - 2k + 2, \lambda_2 = 0, \lambda_3 = 2(k - \lambda - 1), m = 2, n = v$;

(ii) $v_1 = 2v, k_1 = v + 1, \lambda_1 = v - 2(k - \lambda) + 2, \lambda_2 = 2, \lambda_3 = 2(k - \lambda), m = 2, n = v$;

(ii) $v_1 = 3v, k_1 = v, \lambda_1 = v - 2(k - \lambda), \lambda_2 = 0, \lambda_3 = k - \lambda, m = 3, n = v$.

**Proof:** The $srg(v, k, \lambda, \mu)$ is equivalent to 2–class association scheme given by the association matrices $A_1, A_2$ of order $v$. The parameters of the scheme satisfy: $p_{11}^0 = k, p_{22}^0 = v - k - 1, p_{11}^i + p_{22}^i = v - 2k + 2\lambda (i = 1,2); p_{12}^1 = p_{12}^2 = k - \mu$ [See Dey (1986)]. It can be verified that the following are the incidence matrices of STD ($v$) SDRDs:

$(i) N_1 = circ\ (A_1\ A_2)$  $(ii) N_2 = circ\ (I_v + A_1\ I_v + A_2)$  $(iii) N_1 = circ\ (I_v\ A_1\ A_2)$

**Remark 9:** If $v = 2k + 1$, then the SDRD in (i) is $k$–resolvable and that in (ii) is $(k + 1)$–resolvable.

**Example 6:** Using $srg$ (5, 2, 0, 1) in Theorem 8 (ii), we obtain a 3– resolvable SDRD with parameters: $v_1 = 10, k_1 = 6, \lambda_1 = 3, \lambda_2 = 2, \lambda_3 = 4, m = 2, n = 5$.

**Method V: From Affine resolvable SRGDDs**

An $m \times m$ matrix $DS\ (m, s; x)$ with entries from an additive cyclic group $Z_s = \{0, 1, \ldots, s - 1\} (mod\ s); s \geq 2$ is a difference scheme if on a difference $(mod s)$ in any two distinct columns of the matrix each entry of $Z_s$ occurs $x$ times. It can be easily seen that all entries in the first row and first column of a $DS\ (m, s; x)$ can be set 0. The same concept as the difference scheme has been discussed using the terms Difference matrix $D(m, m, s)$ or a generalized Hadamard matrix $GH(s, x)$ by interchanging the roles of rows and columns. Further replacing the group elements of $Z_s$ by corresponding permutation matrices in $DS\ (m, s; x)$, we obtain [See Kadowaki and Kageyama (2009)]:

**Theorem 9:** The existence of a $DS(m, s; x)$ implies the existence of an affine resolvable semi regular GDD with parameters: $v = b = xs^2, r = k = xs, \lambda_1 = 0, \lambda_2 = x, q_1 = 0, q_2 = x, m = xs, n = s(\geq 2)$. (15)

One of the resolutions classes of the GDD with parameters (15) contains blocks as the row of the GD association scheme. Removing this class we obtain:

**Theorem 10:** The existence of a $DS\ (m, s; x)$ implies the existence of an LSR affine resolvable STD ($s$) RD with parameters: $v = xs^2, b = s(xs - 1), r = xs - 1, k = xs, \lambda_1 = 0, \lambda_2 = x - 1, \lambda_3 = x, q_1 = 0, q_2 = x, m = xs, n = s(\geq 2)$. (16)

Further removing $t$ rows of blocks of incidence matrix of the RD with parameters (16) we obtain:

**Corollary 12:** The existence of a $DS\ (m, s; x)$ implies the existence of a resolvable STD ($s$) SR RD with parameters: $v = s(xs - t), b = s(xs - 1), r = xs - 1, k = xs - t, \lambda_1 = 0, \lambda_2 = x - 1, \lambda_3 = x, m = xs - t, n = s(\geq 2)$. (17)

**5. E–optimal RDs**

Bagchi (2004) defined two classes of E–optimal RDs in the class of 1–designs $D$ $(v, b, r, k)$ for $m = 2$:

(1) Balanced of type 1 which satisfies $\lambda_1 = \lambda_2 = \lambda_3 - 1, v > 12$;

(2) Balanced of type 2 which satisfies $\lambda_1 + 1 = \lambda_2 - 1 = \lambda_3, v \geq 10$.

She obtained two series of E–optimal RDs for $v = 2q$ belonging to the two classes, where $q$ is a prime power. Sinha *et al.* (2002a) also obtained a series of E–optimal balanced of type I RDs when $v = 2q, m = 2$. Bose *et al.* (1960) proved that for $n > 6$, there is a pair of MOLS of order $n$. Hence for $m = 2$, Theorem 7 yields:

**Theorem 11:** For $n > 6$ there exist following series of E–optimal balanced of type1 RDs in $D$ $(v, b, r, k)$ with parameters:

(i) $v = 2n, r = n - 1, k = 2, b = n(n-1), \lambda_1 = \lambda_2 = 0, \lambda_3 = 1, m = 2, n$.

(ii) $v = 2n, r = (n-1)^2, k = 2(n-1), b = n(n-1), \lambda_1 = \lambda_2 = 2, \lambda_3 = 3, m = 2, n$ which is the complementary series of (i).

## 6. Parameters of a possible symmetric RD / GDD D when $v = mn$ (even) and $k$ are known

Step I: Obtain $\lambda_i (i = 1,2,3)$ as a solution of the Diophantine equation $(n-1)\lambda_1 + (m-1)\lambda_2 + (n-1)(m-1)\lambda_3 = k(k-1)$, see (1).

Step II: Calculate $\theta_i (i = 1,2,3)$. If a $\theta_i$ is negative, $D$ is not possible. Otherwise we move to Step III.

Step III: (a) If $\theta_1 = 0, \theta_2, \theta_3 > 0$ or $\theta_2 = 0, \theta_1, \theta_3 > 0$, $D$ is SR. (b) If $\theta_3 = 0$, $D$ is singular. (c) If $\theta_1 = \theta_2 = 0; \theta_3 > 0$ then $D$ is LSR RD (d) If $\theta_1, \theta_2, \theta_3 > 0$, then from Theorem 2 (i); a necessary condition for a regular $D$ is: (d1) when $m$ is even and $n$ is odd, $\theta_2$ is a square (d2) when $m$ is odd and $n$ is even, $\theta_1$ is a square, (d3) when $m, n$ are both even, $\theta_1 \theta_2 \theta_3$ is a square.

In the Table 2 given below * over the value of $k$ indicates that there is no SDRD for that value of $k$ as they do not satisfy the above necessary conditions. $[m, n]$ denotes the set of all integers $k$ such that $m \leq k \leq n$. # denotes the set of block size $k$ on which self–dual BIBD/PBIBD (2/3) exists. PBIBD (2/3) stands for partially balanced incomplete block designs with two/three associate classes.

| Table 2: List of Self– dual designs with $k \leq 10$ | | | | | |
|---|---|---|---|---|---|
| No. | v | # | No. | v | # |
| 1 | 8 | [3, 8] | 14 | 28 | [6, 10]– {8} |
| 2 | 9 | [4, 9] | 15 | 30 | [6, 10]– {6, 7} |
| 3 | 10 | [3, 10] | 16 | 33 | [7, 10]–{9} |
| 4 | 14 | [3, 10] | 17 | 35 | [7, 10]–{9} |
| 5 | 15 | [4, 10] | 18 | 36 | [5, 10]–{7, 9} |
| 6 | 16 | [3, 10] | 19 | 38 | {9, 10} |

| 7  | 18  | [6, 10]       | 20 | 40  | {4, 9}          |
|----|-----|---------------|----|-----|-----------------|
| 8  | 20  | [4, 10]–{5*, 6*} | 21 | 42  | [6, 9]– {7*}    |
| 9  | 21  | [5, 10]       | 22 | 44  | [8, 9]          |
| 10 | 22  | [6, 10]–{8, 9} | 23 | 45  | [7, 10]– {8}   |
| 11 | 24  | [5, 10]–{6}   | 24 | 48  | {7, 8}          |
| 12 | 25  | [5, 9]–{6,7}  | 25 | 49  | [7, 10]–{8}    |
| 13 | 27  | [6, 9]        | 26 | 55  | {9, 10}         |
|    |     | –             | 27 | 100 | {9, 10}         |

The following Self dual GDDs/RDs under the range of $r, k \leq 10$ are expected as they satisfy necessary conditions for their existence as described above:

**Table 3: Permissible parameters of GD/RD designs**

| No. | $v$ | $m$ | $n$ | $k$ | $\lambda_1$ | $\lambda_2$ | $\lambda_3$ | Nature of GD/RD |
|-----|-----|-----|-----|-----|-------------|-------------|-------------|------------------|
| 1   | 18  | 3   | 6   | 5   | 2           | 0           | 1           | R RD             |
| 2   | 24  | 3   | 8   | 6   | 2           | 1           | 1           | R GDD            |
| 3   | 26  | 13  | 2   | 9   | 0           | 3           | 3           | R GDD            |
| 4   | 30  | 6   | 5   | 6   | 0           | 2           | 1           | SR RD            |
| 5   | 30  | 3   | 10  | 7   | 2           | 3           | 1           | R RD             |
| 6   | 34  | 2   | 17  | 9   | 2           | 8           | 2           | R GDD            |
| 7   | 36  | 4   | 9   | 9   | 6           | 0           | 1           | SR RD            |
| 8   | 38  | 2   | 19  | 7   | 1           | 6           | 1           | R GDD            |
| 9   | 40  | 5   | 8   | 7   | 2           | 0           | 1           | R RD             |
| 10  | 40  | 5   | 8   | 8   | 4           | 0           | 1           | SR RD            |

## 7. Tables of RDs

In this section two Tables (3 – 4) of RDs in the range of parameters $2 \leq r, k \leq 10$ not found in the Tables of Suen (1989), Sinha *et al.* (1993, 2002*a*) and Kageyama and Sinha (2003) are constructed using the present theorems. Overall efficiency ($E$) of the designs has been calculated using R as given in Kaur *et al.* (2017). $\alpha$ stands for $\alpha -$ resolvability.

**Table 4: Self dual RDs with $2 \leq r, k \leq 10$**

| No. | $v$ | $k$ | $\lambda_1$ | $\lambda_2$ | $\lambda_3$ | $m$ | $n$ | $E$ | $\alpha$ | Nature | Source |
|-----|-----|-----|-------------|-------------|-------------|-----|-----|------|----------|--------|--------|
| 1*  | 6   | 4   | 2           | 2           | 3           | 2   | 3   | 0.89 | 2        | SR     | Cor. 11(ii) |
| 2*  | 10  | 6   | 3           | 2           | 4           | 2   | 5   | 0.92 | 3        | SR     | Th. 8 (ii) |
| 3*  | 12  | 3   | 0           | 0           | 1           | 3   | 4   | 0.68 | 1        | SR     | Cor. 8 |
| 4   | 12  | 8   | 4           | 6           | 5           | 4   | 3   | 0.95 | 2        | SR     | Th. 8 |
| 5*  | 12  | 9   | 6           | 6           | 7           | 4   | 3   | 0.97 | 3        | SR     | Cor. 10 |
| 6   | 14  | 5   | 2           | 2           | 1           | 2   | 7   | 0.85 | -        | R      | Cor. 2 |
| 7*  | 14  | 8   | 4           | 2           | 5           | 2   | 7   | 0.93 | 4        | SR     | Cor. 6 |
| 8   | 14  | 9   | 6           | 6           | 5           | 2   | 7   | 0.95 | -        | R      | Complement of No. 6 |
| 9*  | 18  | 10  | 5           | 2           | 6           | 2   | 9   | 0.95 | 5        | SR     | Th. 8 (ii) |
| 10* | 20  | 4   | 0           | 0           | 1           | 4   | 5   | 0.76 | 1        | SR     | Cor. 8 |
| 11  | 20  | 8   | 2           | 4           | 3           | 4   | 5   | 0.91 | 2        | SR     | Cor. 9 |
| 12* | 21  | 9   | 3           | 3           | 4           | 3   | 7   | 0.93 | 3        | SR     | Example 3 |
| 13  | 22  | 6   | 2           | 0           | 1           | 2   | 11  | 0.86 | -        | R      | Th. 3, $t=3$ |

| | | | | | | | | | | | | |
|---|---|---|---|---|---|---|---|---|---|---|---|---|
| 14 | 22 | 7 | 3 | 2 | 1 | 2 | 11 | 0.88 | - | R | Cor. 2 |
| 15 | 28 | 7 | 2 | 4 | 1 | 4 | 7 | 0.87 | - | R | Cor. 2 |
| 16[*] | 28 | 7 | 0 | 1 | 2 | 7 | 4 | 0.88 | 1 | SR | Cor.12 |
| 17 | 30 | 8 | 3 | 0 | 1 | 2 | 15 | 0.88 | - | R | Th. 3; $t=4$ |
| 18 | 30 | 9 | 4 | 2 | 1 | 2 | 15 | 0.89 | - | R | Cor. 2 |
| 19 | 33 | 8 | 3 | 3 | 1 | 3 | 11 | 0.89 | - | R | Cor. 2 |
| 20 | 35 | 8 | 2 | 5 | 1 | 5 | 7 | 0.88 | - | R | Cor. 2 |
| 21 | 38 | 10 | 4 | 0 | 1 | 2 | 19 | 0.90 | - | R | Th. 3, $t=5$ |
| 22 | 42 | 9 | 2 | 6 | 1 | 6 | 7 | 0.88 | - | R | Cor. 2 |
| 23 | 44 | 9 | 3 | 4 | 1 | 4 | 11 | 0.90 | - | R | Cor. 2 |
| 24[*] | 45 | 9 | 0 | 1 | 2 | 9 | 5 | 0.91 | 1 | SR | Cor. 12 |
| 25 | 45 | 10 | 4 | 3 | 1 | 3 | 15 | 0.90 | - | R | Th. 3, $t=4$ |
| 26 | 49 | 10 | 2 | 7 | 1 | 7 | 7 | 0.88 | - | R | Cor. 2 |
| 27 | 55 | 10 | 3 | 5 | 1 | 5 | 11 | 0.90 | - | R | Cor. 2 |
| 28 | 60 | 10 | 3 | 2 | 1 | 4 | 15 | 0.90 | - | R | Th. 3, $t=4$ |

Table 5: Table of RDs with $\lambda_3 > \lambda_1, \lambda_3 > \lambda_2$ and $2 \leq r, k \leq 10$

| No. | $v$ | $b$ | $r$ | $k$ | $\lambda_1$ | $\lambda_2$ | $\lambda_3$ | $m$ | $n$ | $E$ | $\alpha$ | Nature | Source |
|---|---|---|---|---|---|---|---|---|---|---|---|---|---|
| 1 | 6 | 14 | 7 | 3 | 0 | 3 | 4 | 3 | 2 | 0.76 | 1 | SR | Cor. 12 |
| 2 | 8 | 28 | 7 | 2 | 0 | 1 | 2 | 2 | 4 | 0.53 | 1 | SR | Cor. 12 |
| 3 | 8 | 14 | 7 | 4 | 0 | 3 | 4 | 4 | 2 | 0.84 | 1 | SR | Cor. 12 |
| 4 | 10 | 45 | 9 | 2 | 0 | 1 | 2 | 2 | 5 | 0.52 | 1 | SR | Cor. 12 |
| 5 | 10 | 14 | 7 | 5 | 0 | 3 | 4 | 5 | 2 | 0.88 | 1 | SR | Cor. 12; DS (8, 2; 4) |
| 6 | 12 | 28 | 7 | 3 | 0 | 1 | 2 | 3 | 4 | 0.70 | 1 | SR | Cor. 12; DS (8, 4; 2) |
| 7 | 15 | 45 | 9 | 3 | 0 | 1 | 2 | 3 | 5 | 0.70 | 1 | SR | Cor. 12 |
| 8 | 20 | 45 | 9 | 4 | 0 | 1 | 2 | 4 | 5 | 0.78 | 1 | SR | Cor. 12 |
| 9 | 22 | 55 | 10 | 4 | 0 | 1 | 2 | 11 | 2 | 0.78 | 2 | SR | Example 4 |
| 10 | 25 | 45 | 9 | 5 | 0 | 1 | 2 | 5 | 5 | 0.83 | 1 | SR | Cor. 12 |
| 11 | 30 | 45 | 9 | 6 | 0 | 1 | 2 | 6 | 5 | 0.86 | 1 | SR | Cor. 12 |
| 12 | 35 | 45 | 9 | 7 | 0 | 1 | 2 | 7 | 5 | 0.88 | 1 | SR | Cor. 12 |
| 13 | 40 | 45 | 9 | 8 | 0 | 1 | 2 | 8 | 5 | 0.89 | 1 | SR | Cor. 12 |
| 14 | 50 | 45 | 9 | 10 | 0 | 1 | 2 | 10 | 5 | 0.92 | 1 | LSR | Th. 10; DS (10,5; 2) |

Design numbers 2 and 10 of Table 4 are constructed using srg (5, 2, 0, 1) and srg (9, 4, 1, 2) respectively. A Table of $srg(v, k, \lambda, \mu)$ may be found in Brouwer (2007). The design no. 14 of Table 5 is affine resolvable. Some construction methods of difference schemes (generalized Hadamard matrices) may be found in Dey and Mukerjee (1999) and Hedayat *et al.* (1999).

## 8. Conclusion

In this paper we have classified RDs into four classes and constructed some of them. Two classes namely LR and LSR have been identified which have received little attention of authors, one being suitable for confounding. Ionin and Shrikhande (2006) addressed the global decomposition of an SBIBD $D$ which is a tactical decomposition of the incidence matrix $N$ of $D$ as $N =$

$[N_{ij}]_{\substack{1 \leq i \leq s \\ 1 \leq j \leq t}}$ where each $N_{ij}$ is the incidence matrix of some SBIBDs. SDRD has some analogy with SBIBD. Here in every construction we obtain a matrix $N$ whose blocks are square $(0,1)-$ matrices such that $N$ becomes the incidence matrix of a RD. The method is the reverse of the well–known tactical decomposition of the incidence matrix of a known design. Tactical decomposable designs are of great interest because of their connections with automorphisms of designs, see Bekar *et al.* (1982). Saurabh and Singh (2020), Saurabh *et al.* (2021), Saurabh and Sinha (2021, 2022) have constructed some series of GDDs and Latin square designs using this approach. Some new series of RDs with global decomposable and self–dual properties are obtained. Earlier authors have also constructed RDs whose incidence matrices contain incidence matrices of SBIBDs. These methods are scattered over the literature. We have unified and generalized these approaches and constructed some more series of RDs. A series of E – optimal RDs is also constructed, see Theorem 11.

The RDs are useful as factorial experiments, having factorial balance as well as orthogonality [See Dey (1986), Section 6.5.3]. The RDs with $\lambda_3$ bigger than $\lambda_1$ and $\lambda_2$ when used as the $m \times n$ complete confounded factorial experiments; the loss of information on the main effect is small. Such designs are suitable for confounding experiments [See Suen (1989)]. The RDs with asterisk marks in Table 4 have $\lambda_3 > \lambda_1, \lambda_2$ and all the RDs in Table 5 have $\lambda_3 > \lambda_1, \lambda_2$. The RDs in Corollaries 4 and 5 are suitable for confounding if $m = 2$, RDs in Corollary 6 are suitable if $m = 2,3; t > 1$ whereas RDs in Corollary 7 are suitable for $\forall t$ when $m = 2,3$. If $D$ is LSR RD then $\lambda_3 > \lambda_1, \lambda_2$, see Proposition 2. Theorem 10 and Corollaries 8,10,12 are series of LSR RDs. Sinha *et al.* (1993) [See Theorem 2.5] also obtained a series of LSR RDs. Their design numbers 13, 19, 39, 41, 56, 59, 60, 63 and 65 given in Table 1 are LSR RDs. We have obtained only one series of LSR RDs and only one series of LR RDs. The series of LR RDs yields only one design under the range $2 \leq r, k \leq 10$. It would be interesting and challenging to construct new series of LR, LSR RDs and regular E–optimal RDs on even $v$ using $(0, 1)$ – linear combination of association matrices.

The rectangular designs have also been used in the construction of balanced and orthogonal arrays, see Sinha *et al.* (2002*b*). Sinha *et al.* (2002*a*) stated the Theorem 2.5 for the construction of RD from BIBD. Unfortunately RDs numbered 12, 15, 20, 21, 27, 28, 32, 33 of the Table 1 in Sinha *et al.* (2002) constructed using Theorem 2.5 do not exist since the eigenvalues of their incidence matrices are negative. Still there are several values of even $v$ and $k$ on which there is no self–dual BIBD/PBIBD (2/3), see Table 2.